\documentclass[preprint,english]{siamltex}
\usepackage[T1]{fontenc}
\usepackage[latin9]{inputenc}
\usepackage{float}
\usepackage{amsmath}
\usepackage{amssymb}
\usepackage{graphicx}
\usepackage{xcolor}

\makeatletter
\newtheorem{thm}{\protect\theoremname}
\newtheorem{rem}[thm]{\protect\remarkname}

\def\la{\langle}
\def\ra{\rangle}

\usepackage{fullpage}

\makeatother

\usepackage{babel}
\providecommand{\remarkname}{Remark}
\providecommand{\theoremname}{Theorem}

\title{A quasi-conservative dynamical low-rank algorithm for the Vlasov equation}

\author{Lukas Einkemmer\footnotemark[1]\ \footnotemark[3] \and Christian Lubich\footnotemark[1]}

\begin{document}

\maketitle

\renewcommand{\thefootnote}{\fnsymbol{footnote}}
\footnotetext[1]{Mathematisches Institut,
       Universit\"at T\"ubingen,
       Auf der Morgenstelle 10,
       D--72076 T\"ubingen,
       Germany. Email: {\tt \{einkemmer,lubich\}@na.uni-tuebingen.de}}
 \footnotetext[3]{Department of Mathematics, University of Innsbruck, Austria}
\renewcommand{\thefootnote}{\arabic{footnote}}

\begin{abstract} Numerical methods that approximate the solution of the Vlasov--Poisson equation by a low-rank representation have been considered recently. These methods can be extremely effective from a computational point of view, but contrary to most Eulerian Vlasov solvers, they do not conserve mass and momentum, neither globally nor in respecting the corresponding local conservation laws. This can be a significant limitation for intermediate and long time integration. In this paper we propose a numerical algorithm that overcomes some of these difficulties and demonstrate its utility by presenting numerical simulations.
\end{abstract}

\section{Introduction}

Many plasma systems that are of interest in applications (such as in
magnetic confined fusion or astrophysics) cannot be adequately described
by fluid models. Instead kinetic models have to be employed. Since
these models are posed in a $2d$-dimensional ($d=1,2,3$) phase
space, numerically solving kinetic equations on a grid is extremely
expensive from a computational point of view. Thus, traditionally,
particle methods have been employed extensively to approximate these
types of problems (see, for example, \cite{verboncoeur2005particle}).
However, particle methods suffer from excessive noise that makes it,
for example, difficult to resolve regions with low phase space density.
Due to the increase in computer performance, methods that directly
discretize phase space, the so-called Eulerian approach, have recently
seen increased interest \cite{sonnendrucker1999,filbet2003comparison,crouseilles2009forward,sircombe2009valis,qiu2011,rossmanith2011,crouseilles2011,einkemmer2014,crouseilles2015hamiltonian,einkemmer2015,crouseilles2016asymptotic,einkemmer2017study}.
However, performing these simulations in higher dimensions is still
extremely expensive. As a consequence, much effort has been devoted
to efficiently implement these methods on high performance computing
systems \cite{rozar2013,bigot2013,einkemmer2015,latu2007gyrokinetic,mehrenberger2013vlasov,einkemmer2016mixed,crouseilles2009parallel,einkemmer2018comparison}.

More recently, methods that use a low-rank approximation have emerged.
In \cite{Ehrlacher2017,Kormanna} the Vlasov equation is first discretized in time and/or space, and then low-rank algorithms are applied to the discretized system.
A different approach is taken
in \cite{einkemmer2018low}, where a low-rank projector-splitting is on top of the procedure. That is, the low-rank algorithm is applied before any time or space discretization is chosen. This results in small systems of $d$-dimensional advection equations (in either the space or the velocity variables, in an alternating fashion) that are then solved by spectral or semi-Lagrangian methods. The advantage
of this approach is that the evolution equations that need to be solved
numerically are directly posed in terms of the degrees of freedom
of the low-rank representation. Thus, no intermediate tensors have
to be constructed and no tensor truncation algorithms have
to be employed. This also leads to increased flexibility in the choice
of the time and space discretization methods.

Computing numerical solutions of high-dimensional evolutionary
partial differential equations by {\it dynamical low-rank approximation}
has only recently been considered for kinetic problems \cite{einkemmer2018low,einkemmer2018weakly}.
However, such algorithms have been investigated extensively in quantum
mechanics; see, in particular, \cite{meyer90tmc,meyer09mqd} for the
MCTDH approach to molecular quantum dynamics in the chemical physics
literature and \cite{Lubich2008,lubich15tii} for a computational
mathematics point of view of this approach. Some uses of dynamical
low-rank approximation in areas outside quantum mechanics are described
in \cite{Nonnenmacher2008,jahnke2008dynamical,Mena2017,Musharbash2018}.
In a general mathematical setting, dynamical low-rank approximation
has been studied in \cite{Koch2007,Koch2010,lubich13dab}. A major
algorithmic advance for the time integration was achieved with the
projector-splitting methods first proposed in \cite{Lubich2014} for
matrix differential equations and then developed further for various
tensor formats in \cite{lubich15tii,lubich15tio,Haegeman2016,Kieri2016,Lubich2018}.
In contrast to standard time-stepping methods, the projector-splitting methods have been shown to be robust to the typical presence of small singular values in the low-rank approximation~\cite{Kieri2016}.
The approach in \cite{einkemmer2018low,einkemmer2018weakly} and in the present paper is
based on an adaptation of the projector-splitting method of \cite{Lubich2014} to kinetic equations. 

While low-rank approximations can be very effective from a computational
point of view, they destroy much of the physical structure of the
problem under consideration. 
Important physical invariants, such as mass and momentum, are no longer conserved. Perhaps even more problematic is
that the low-rank approximation does not take the corresponding local
conservation laws into account. This can be a significant issue if
these algorithms are to be used for long or even intermediate time
integration. 

This situation is in stark contrast with the state of the art for
Eulerian Vlasov solvers, where significant research has been conducted
to conserve certain physical properties of the exact solution \cite{filbet2003comparison,sircombe2009valis,qiu2010conservative,crouseilles2010conservative,crouseilles2015hamiltonian,crouseilles2016asymptotic,einkemmer2017study}.
In particular, methods that conserve mass and momentum are commonly
employed. However, to the best of our knowledge, no low-rank algorithms
are available that are able to conserve even linear invariants.
Furthermore, it has recently been proposed to use low-rank
numerical methods to solve fluid problems \cite{einkemmer2018weakly}.
Also in this setting conservation of mass and momentum, a hallmark of
traditional fluid solvers, is, of course, of great interest.

In this paper we will consider the Vlasov\textendash Poisson equation

\begin{equation}\label{eq:vlasov-poisson}
\begin{aligned} 
& \partial_{t}f(t,x,v)+v\cdot\nabla_{x}f(t,x,v)-E(f)(x)\cdot\nabla_{v}f(t,x,v)=0 \\
& \nabla\cdot E(f)(x)=-\int f(t,x,v)\,\mathrm{d}v+1,\qquad\;\;\nabla\times E(f)(x)=0,
\end{aligned} \end{equation}which models the time evolution of a collisionless plasma in the electrostatic
regime. This equation has an infinite number of invariants (Casimir
invariants). Here we will consider the linear invariants of mass and
momentum and the corresponding local conservation laws. In section
\ref{sec:Dynamical-low-rank-splitting} we will introduce the necessary
notation and describe the dynamical low-rank splitting algorithm for
the Vlasov equation that was proposed in \cite{einkemmer2018low}.
We then derive a modification of that numerical method such that a
projected version of the continuity and momentum balance equation
is satisfied (section \ref{subsec:local-conservation}). Subsequently
we will discuss the global conservation of mass and momentum in section
\ref{subsec:Global-conservation}. We will then consider the efficient
implementation of these methods (section \ref{sec:Efficient-implementation}).
Finally, in section \ref{sec:numerical-results} we present numerical
results for the Vlasov\textendash Poisson equation. In particular,
we will demonstrate the efficiency of the proposed algorithms for
a two-stream instability.

\section{A low-rank projector-splitting integrator\label{sec:Dynamical-low-rank-splitting}}

We will start by summarizing the low-rank projector splitting integrator
for the Vlasov\textendash Poisson equation introduced in \cite{einkemmer2018low}.
It should be duly noted that this algorithm neither respects the local conservation
laws associated with mass or momentum, nor conserves
mass or momentum globally (this is also true for low-rank algorithms in \cite{Ehrlacher2017,Kormanna}).

We seek an approximation to the Vlasov\textendash Poisson equation
(\ref{eq:vlasov-poisson}) in the following form:
\[
f(t,x,v)=\sum_{i,j=1}^{r}X_{i}(t,x)S_{ij}(t)V_{j}(t,v),
\]
with real coefficients $S_{ij}(t)\in\mathbb{R}$ and with functions $X_{i}(t,x)$ and $V_{j}(t,v)$ that are orthonormal:
$$
\langle X_{i},X_{k}\rangle_x=\delta_{ik} \quad\text{ and }\quad \langle V_{j},V_{l}\rangle_v=\delta_{jl},
$$
where $\la\cdot,\cdot\ra_x$ and $\la\cdot,\cdot\ra_v$ are the inner products on $L^{2}(\Omega_{x})$
and $L^{2}(\Omega_{v})$, respectively.
The dependence of $f$ on the phase
space variables $(x,v)\in\Omega=\Omega_{x}\times\Omega_{v}\subset\mathbb{R}^{2d}$
is approximated by the functions $\left\{ X_{i}\colon i=1,\ldots,r\right\} $
and $\left\{ V_{j}\colon j=1,\ldots,r\right\} $, which depend only
on the separated variables $x\in\Omega_{x}$ and $v\in\Omega_{v}\subset\mathbb{R}^{d}$ ,
respectively. Such an approach is efficient if the rank $r$ can be
chosen much smaller compared to the number of grid points used to
discretize $X_{i}$ and $V_{j}$ in space.

The dynamics of the Vlasov\textendash Poisson equation is constrained
to the corresponding low-rank manifold by replacing  \eqref{eq:vlasov-poisson} with an evolution equation
\[
\partial_{t}f=-P(f)\left(v\cdot\nabla_{x}f-E(f)\cdot\nabla_{v}f\right),
\]
where $P(f)$ is the orthogonal projector onto the manifold. The projector
can be written as 
\begin{equation}
P(f)g=P_{\overline{V}}g-P_{\overline{V}}P_{\overline{X}}g+P_{\overline{X}}g,\label{eq:projector-P(f)}
\end{equation}
where $P_{\overline{X}}$ is the orthogonal projector onto the vector
space $\overline{X}=\text{span}\left\{ X_{i}\colon i=1,\ldots,r\right\} $
and $P_{\overline{V}}$ is the orthogonal projector onto the vector
space $\overline V=\text{span}\left\{ V_{j}\colon j=1,\ldots,r\right\} $. Then,
as first suggested in \cite{Lubich2014}, the dynamics is split into
the three terms of equation (\ref{eq:projector-P(f)}). In the simplest
case, the first-order Lie-Trotter splitting, we solve the equations
\begin{align}
\partial_{t}f & =-P_{\overline{V}}(v\cdot\nabla_{x}f-E(f)\cdot\nabla_{v}f)\label{eq:splitting-eq1}\\
\partial_{t}f & =+P_{\overline{V}}P_{\overline{X}}(v\cdot\nabla_{x}f-E(f)\cdot\nabla_{v}f)\label{eq:splitting-eq2}\\
\partial_{t}f & =-P_{\overline{X}}(v\cdot\nabla_{x}f-E(f)\cdot\nabla_{v}f)\label{eq:splitting-eq3}
\end{align}
one after the other. Now, let us define
\[
K_{j}(t,x)=\sum_{i}X_{i}(t,x)S_{ij}(t),\qquad\qquad L_{i}(t,v)=\sum_{j}S_{ij}(t)V_{j}(t,v).
\]
The advantage of the splitting scheme then becomes that equation (\ref{eq:splitting-eq1})
only updates $K_{j}$ (the $V_{j}$ stay constant during that step),
equation (\ref{eq:splitting-eq2}) only updates $S_{ij}$ (the $X_{i}$
and $V_{j}$ stay constant during that step), and equation (\ref{eq:splitting-eq3})
only updates $L_{i}$ (the $X_{i}$ stay constant during that step).
The corresponding evolution equations are derived in \cite{einkemmer2018low}
and are of the following form:
\begin{align}
\partial_{t}K_{j}(t,x) & =-\sum_{l}c_{jl}^{1}\cdot\nabla_{x}K_{l}(t,x)+\sum_{l}c_{jl}^{2}\cdot E(K)(t,x)K_{l}(t,x)\label{eq:evol-K}\\
\partial_{t}S_{ij}(t) & =\sum_{k,l}(c_{jl}^{1}\cdot d_{ik}^{2}-c_{jl}^{2}\cdot d_{ik}^{1}[E(S(t))])S_{kl}(t)\label{eq:evol-S}\\
\partial_{t}L_{i}(t,v) & =\sum_{k}d_{ik}^{1}[E(L(t,\cdot))]\cdot\nabla_{v}L_{k}(t,v)-\sum_{k}(d_{ik}^{2}\cdot v)L_{k}(t,v).\label{eq:evol-L}
\end{align}
The coefficients $c_{jl}^{1}$, $c_{jl}^{2}$ and $d_{ik}^{1}$, $d_{ik}^{2}$ are vector-valued
but constant in  space and, with the exception of $d_{ik}^1$, also constant in time. They are given by integrals over $\Omega_v$ and $\Omega_x$, respectively; see \cite[Section~2]{einkemmer2018low} for the details. 

Assuming that the initial value is
represented as $f^{0}(x,v)=\sum_{i,j}X_{i}^{0}(x)S_{ij}^{0}V_{j}^{0}(v)$,
the algorithm with time step size $\tau$ then proceeds in the following
three steps.

\textbf{Step 1: }Solve equation (\ref{eq:evol-K}) with initial value
$K_{j}(0,x)=K_{j}^{0}=\sum_{i}X_{i}^{0}(x) S_{ij}^{0}.$ Then perform
a QR decomposition of $K^{1}=[K_{1}(\tau,\cdot),\dots,K_{r}(\tau,\cdot)]$
to obtain $X_{i}^{1}$ and $\widehat S_{ij}^{1}$.

\textbf{Step 2: }Solve equation (\ref{eq:evol-S}) with initial value
$S_{ij}(0)=\widehat S_{ij}^{1}$ to obtain $\widetilde S_{ij}^{0}=S_{ij}(\tau)$.

\textbf{Step 3: }Solve equation (\ref{eq:evol-L}) with initial value
$L_{i}(0,v)=L_{i}^{0}=\sum_{j}\widetilde S_{ij}^{0}V_{j}^{0}$. Then perform
a QR decomposition of $L^{1}=[L_{1}(\tau,\cdot),\dots,L_{r}(\tau,\cdot)]$
to obtain  $V_{j}^{1}$ and $S_{ij}^{1}$.

The output of the algorithm is then the low-rank representation 
\[
f(\tau,x,v)\approx f^{1}(x,v)=\sum_{i,j}X_{i}^{1}S_{ij}^{1}V_{i}^{1}.
\]
For a detailed derivation of this algorithm the reader is referred
to \cite{einkemmer2018low}. We note that the extension to second
order Strang splitting is immediate.

\section{Local conservation\label{subsec:local-conservation}}

The Vlasov\textendash Poisson equation (\ref{eq:vlasov-poisson})
satisfies the continuity equation
\begin{equation}
\partial_{t}\rho(t,x)+\nabla\cdot(\rho(t,x)u(t,x))=0\label{eq:continuity-eq}
\end{equation}
and the momentum balance equation
\begin{equation}
\partial_{t}\bigl(\rho(t,x)u(t,x)\bigr)+\nabla\cdot(\rho(t,x)u(t,x)\otimes u(t,x))=-E(t,x)\rho(t,x),\label{eq:momentum-balance-eq}
\end{equation}
where
\[
\rho(t,x)=\int f(t,x,v)\,\mathrm{d}v,\qquad\rho(t,x)u(t,x)=\int vf(t,x,v)\,\mathrm{d}v.
\]
From these equations, global conservation of mass and momentum is
easily obtained by integrating in $x$. Without the projection operators,
equations (\ref{eq:splitting-eq1})--(\ref{eq:splitting-eq3})
would satisfy the continuity equation (\ref{eq:continuity-eq}) and
the momentum balance equation (\ref{eq:momentum-balance-eq}). Overall
this would ensure that the splitting scheme (without projection operators)
respects the local conservation laws for mass and momentum. However,
it can easily be seen that the projection operators destroy this property.
In addition, as has already been pointed out in \cite{einkemmer2018low},
global conservation of mass and momentum is lost as well.

A crucial observation that enables the following numerical method
is the observation that the conserved quantities only depend on $x$.
While we cannot modify the algorithm such that the conservation laws are satisfied exactly (while keeping $V_j$ constant in Step 1 and $X_i$ constant in Step 3, and both $X_i$ and $V_j$ constant in Step 2), our goal is to derive
a numerical method that satisfies the {\it projected} conservation laws for
mass and momentum
\begin{equation} \label{proj-cons}
P_{\overline{X}}\left(\partial_{t}\rho+\nabla\cdot(\rho u)\right)=0,\qquad\qquad P_{\overline{X}}(\partial_{t}\rho+\nabla\cdot(\rho u\otimes u)+E\rho)=0.
\end{equation}
The idea is to add to (\ref{eq:splitting-eq1})--(\ref{eq:splitting-eq3}) corrections of the form
\begin{equation}
\sum_{i,j}\lambda_{ij}X_{i}V_{j},\label{eq:ansatz}
\end{equation}
where the coefficients $\lambda_{ij}$ are determined such that the projected continuity equation
and the projected momentum balance hold true. This results in an overdetermined
system for the $\lambda_{ij}$ for which we seek the smallest solution in the Euclidean norm.

One might object at this point and argue that such a correction is
unnecessarily restrictive. Certainly, one could envisage that for equation
(\ref{eq:evol-K}) and (\ref{eq:evol-L}) an arbitrary function of
$x$ and $v$, respectively, could be used as the correction. Unfortunately,
as we will describe in more detail in Remark \ref{rem:Rcorrection},
this would introduce, for example, non-zero values in the density
function at high velocities. This, clearly unphysical, artefact then
pollutes the numerical solution. Thus, the benefit of the ansatz given
in equation (\ref{eq:ansatz}) is that the $X_{i}$ and $V_{j}$,
which are already used to represent the numerical solution, are also
used for the correction. 
Since the
algorithm adapts the functions $X_{i}$ and $V_{j}$ in accordance
with the solution, the artefact described above is avoided. This
behavior is confirmed by numerical simulation.

In the following, the correction given in (\ref{eq:ansatz}) will
be made precise for the three steps of the splitting algorithm.

\textbf{Step 1: }We replace the evolution equation (\ref{eq:evol-K})
by
\begin{align}
\partial_{t}K_{j} & =\left\langle V_{j}^{0},F(f)+\sum_{k,l}\lambda_{kl}X_{k}^{0}V_{l}^{0}\right\rangle _{v}\nonumber \\
 & =\langle V_{j}^{0},F(f)\rangle_{v}+\sum_{k}\lambda_{kj}X_{k}^{0}\label{eq:evolution-corr-K}
\end{align}
with $F(f)=-v\cdot\nabla_{x}f+E(f)\cdot\nabla_{v}f$ for $f(t,x,v)=\sum_l K_l(t,x) V_l^0(v)$, and $\lambda_{kl}$ is
yet to be determined. Note that the $V_{l}$ are constant during
that time step, and hence $F(f)$ only depends on the $K_l$. Now,
we impose
\begin{equation}
0=P_{\overline{X}^0}(\partial_{t}\rho+\nabla\cdot(\rho u))=\sum_{i}X_{i}^{0}\left[\sum_{j}\lambda_{ij}\alpha_{j}+\sum_{j}\langle X_{i}^{0}V_{j}^{0},F(f)\rangle_{xv}\alpha_{j}+\langle X_{i}^{0},\nabla\cdot(\rho u)\rangle_{x}\right],\label{eq:continuity-K}
\end{equation}
where $\alpha_{j}=\int V_{j}^{0}\,\mathrm{d}v$, and
\begin{align}
0&=P_{\overline{X}^0}\left(\partial_{t}(\rho u)+\nabla\cdot(\rho u\otimes u)+E\rho\right)
\nonumber
\\&=\sum_{i}X_{i}^{0}\left[\sum_{j}\lambda_{ij}\beta_{j}+\sum_{j}\langle X_{i}^{0}V_{j}^{0},F(f)\rangle_{xv}\beta_{j}+\langle X_{i}^{0},\nabla\cdot(\rho u\otimes u)\rangle_{x}+\langle X_{i}^{0},E(f)\rho\rangle_{x}\right],\label{eq:momentumbal-K}
\end{align}
where $\beta_{j}=\int vV_{j}^{0}\,\mathrm{d}v\in\mathbb{R}^{d}$.
Together, equations (\ref{eq:continuity-K}) and (\ref{eq:momentumbal-K})
yield $(1+d)r$ linear equations for the $r^{2}$ unknowns $\lambda_{ij}$ (We suppose $r\ge 1+d$ in the following).
Since the equations for different $i$ decouple, this allows us to
put this into matrix form as follows, with the row vector $\alpha=(\alpha_1,\dots,\alpha_r)$ and with the $d\times r$ matrix 
$\beta=(\beta_1,\dots,\beta_r)$:
\begin{equation}
\left[\begin{array}{c}
\alpha\\
\beta
\end{array}\right]\lambda_{i(\cdot)}=\left[\begin{array}{c}
b_{i}\\
d_{i}
\end{array}\right]\label{eq:linear-system}
\end{equation}
with
\begin{align*}
b_{i} & =-\sum_{j}\langle X_{i}^{0}V_{j}^{0},F(f)\rangle_{xv}\alpha_{j}-\langle X_{i}^{0},\nabla\cdot(\rho u)\rangle_{x}\\
d_{i} & =-\sum_{j}\langle X_{i}^{0}V_{j}^{0},F(f)\rangle_{xv}\beta_{j}-\langle X_{i}^{0},\nabla\cdot(\rho u\otimes u)\rangle_{x}-\langle X_{i}^{0},E(f)\rho\rangle_{x}.
\end{align*}
These systems of equations have (multiple) solutions if the rows of the matrix $[\alpha;\beta]$ are linearly independent. In order to minimize the magnitude
of the correction that is applied, we seek the solution with the smallest
Euclidean norm. This can be done easily and at negligible cost as the matrix
is only of size $(1+d)\times r$.

It is still necessary to compute the right hand side. We have
\[
\nabla\cdot(\rho u)=\sum_{j}\nabla K_{j}\cdot\beta_{j},\qquad\nabla\cdot(\rho u\otimes u)=\sum_{j}\nabla K_{j}\cdot\gamma_{j},
\]
where $\gamma_{j}=\int(v\otimes v)V_{j}^{0}\,\mathrm{d}v$. Since
$E$ and $\rho$ have to be computed in any case and $\alpha$, $\beta$,
$\gamma$, on modern computer architectures, can be computed alongside
the coefficients $c^{1}$ and $c^{2}$ at (almost) no extra cost,
only the projections in $x$ are of any concern from a computational
point of view. These require the computation of $r$ integrals and
consequently $\mathcal{O}\left(rn^{d}\right)$ arithmetic
operations when $n$ quadrature points are used in each coordinate direction. 

\textbf{Step 2: }We replace the evolution equation (\ref{eq:evol-S})
by
\begin{align}
\partial_{t}S_{ij} & =-\left\langle X_{i}^{1}V_{j}^{0},F(f)+\sum_{k,l}\lambda_{kl}X_{k}^{1}V_{l}^{0}\right\rangle _{xv}\nonumber \\
 & =-\langle X_{i}^{1}V_{j}^{0},F(f)\rangle_{xv}-\lambda_{ij}\label{eq:evolution-corr-S}
\end{align}
for $f(t,x,v)=\sum_{k,l} X_k^0(x) \,S_{kl}(t)\, V_l^0(v)$, so that $F(f)$ depends only on the $S_{kl}$, and
where $\lambda_{ij}$ is yet to be determined. Then we impose the constraints
\begin{equation}
0=P_{\overline{X}^1}\left(\partial_{t}\rho-\nabla\cdot(\rho u)\right)=-\sum_{i}X_{i}^{1}\left[\sum_{j}\lambda_{ij}\alpha_{j}+\sum_{j}\langle X_{i}^{1}V_{j}^{0},F(f)\rangle_{xv}\alpha_{j}+\langle X_{i}^{1},\nabla\cdot(\rho u)\rangle_{x}\right]\label{eq:continuity-S}
\end{equation}
with $\alpha_{j}=\int V_{j}^{0}\,\mathrm{d}v$ and
\begin{align}
0&=P_{\overline{X}^1}\left(\partial_{t}(\rho u)-\nabla\cdot\left(\rho u\otimes u\right)-E\rho\right) \nonumber
\\
&=-\sum_{i}X_{i}^{1}\left[\sum_{j}\lambda_{ij}\beta_{j}+\sum_{j}\langle X_{i}^{1}V_{j}^{0},F(f)\rangle_{xv}\beta_{j}+\langle X_{i}^{1},\nabla\cdot(\rho u\otimes u)\rangle_{x}+\langle X_{i}^{1},E\rho\rangle_{x}\right]\label{eq:momentumbal-S}
\end{align}
with $\beta_{j}=\int vV_{j}^0\,\mathrm{d}v$. Equations (\ref{eq:continuity-S})
and (\ref{eq:momentumbal-S}) yield $(1+d)r$ linear equations for
the $r^{2}$ unknowns $\lambda_{ij}$. Since the equations for different
$i$ decouple, we once again can put this into the form given by equation
(\ref{eq:linear-system}). The only difference lies in the right hand
side which is computed as follows:

\begin{align*}
b_{i} & =-\sum_{j}\langle X_{i}^{1}V_{j}^{0},F(f)\rangle_{xv}\alpha_{j}-\langle X_{i}^{1},\nabla\cdot(\rho u)\rangle_{x}\\
d_{i} & =-\sum_{j}\langle X_{i}^{1}V_{j}^{0},F(f)\rangle_{xv}\beta_{j}-\langle X_{i}^{1},\nabla\cdot(\rho u\otimes u)\rangle_{x}-\langle X_{i}^{1},E\rho\rangle_{x}.
\end{align*}
As before, we seek the solution that minimizes the Euclidean norm of
the $\lambda_{ij}$. This can be done efficiently as we only have to solve $r$
systems of size $(1+d)\times r$. Computing the right-hand side requires
$\langle X_{i}^{1}V_{j}^{0},F(f)\rangle_{xv}$, which has to be computed
to conduct this splitting step in any case. Thus, only the projections
in $x$ remain. As noted above, they can be computed in $\mathcal{O}\left(rn^{d}\right)$
arithmetic operations when $n$ quadrature points are used in each coordinate direction.

\textbf{Step 3: }We replace the evolution equation (\ref{eq:evol-K})
by
\begin{align}
\partial_{t}L_{i} & =\left\langle X_{i}^{1},F(f)+\sum_{kl}\lambda_{kl}X_{k}^{1}V_{l}^{0}\right\rangle _{x}\nonumber \\
 & =\langle X_{i}^{1},F(f)\rangle_{x}+\sum_{l}\lambda_{il}V_{l}^{0}
 \label{eq:evolution-corr-L}
\end{align}
for $f(t,x,v)=\sum_{k,l} X_k^1(x) \,L_k(t,v)$, so that $F(f)$ depends only on the functions $L_k$, and
where $\lambda_{ij}$ is yet to be determined. Then we impose the constraints

\begin{equation}
0=P_{\overline{X}^1}\left(\partial_{t}\rho+\nabla\cdot(\rho u)\right)=\sum_{i}X_{i}^{1}\left[\sum_{l}\lambda_{il}\alpha_{l}+\langle X_{i}^{1},F(f)\rangle_{xv}+\langle X_{i}^{1},\nabla\cdot(\rho u)\rangle_{x}\right]\label{eq:continuity-L}
\end{equation}
with $\alpha_{l}=\int V_{l}^{0}\,\mathrm{d}v$ and
\begin{align}\nonumber
0&=P_{\overline{X}^1}\left(\partial_{t}(\rho u)+\nabla\cdot\left(\rho u\otimes u\right)+E\rho\right)
\\
&=\sum_{i}X_{i}^{1}\left[\sum_{l}\lambda_{il}\beta_{l}+\langle X_{i}^{1},vF(f)\rangle_{xv}+\langle X_{i}^{1},\nabla\cdot(\rho u\otimes u)\rangle_{x}-\langle X_{i}^{1},E\rho\rangle_{x}\right]\label{eq:momentumbal-L}
\end{align}
with $\beta_{l}=\int vV_{l}^{0}\,\mathrm{d}v$. As before, equations
(\ref{eq:continuity-L}) and (\ref{eq:momentumbal-L}) yield $(1+d)r$
linear equations for the $r^{2}$ unknowns $\lambda_{ij}$. We can once
again put this into the form given by equation (\ref{eq:linear-system})
with right-hand side
\begin{align*}
b_{i} & =-\langle X_{i}^{1},F(f)\rangle_{xv}-\langle X_{i}^{1},\nabla\cdot(\rho u)\rangle_{x}\\
d_{i} & =-\langle X_{i}^{1},vF(f)\rangle_{xv}-\langle X_{i}^{1},\nabla\cdot(\rho u\otimes u)\rangle_{x}+\langle X_{i}^{1},E\rho\rangle_{x}.
\end{align*}
As before, this can be done efficiently as the matrix involved is
small and the right-hand side can be efficiently computed alongside
the coefficients that are needed for the low-rank splitting algorithm.

Note that in the third step we have
\begin{align*}
b_{i} & =\sum_{k}\langle X_{i}^{1},\nabla X_{k}^{1}\rangle_{x}\cdot\int vL_{k}\,\mathrm{d}v-\sum_{k}\langle X_{i}^{1},EX_{k}^{1}\rangle_{x}\int\nabla_{v}L_{k}\,\mathrm{d}v-\sum_{k}\langle X_{i}^{1},\nabla X_{k}^{1}\rangle_x\cdot\int vL_{k}\,\mathrm{d}v\\
 & =0,
\end{align*}
where we have assumed that the $L_{k}$ go to zero as $\vert v\vert\to\infty$. Thus,
step 3 already satisfies the continuity equation.

\begin{rem}
At first sight it looks more natural to use $K_{j}$ and $L_{i}$
instead of $X_{j}^{0}$ and $V_{i}^{0}$ in step 1 and 3. These are
the quantities that are updated in that step of the algorithm. The
correction would then also reflect the corresponding changes that
occur as the subflows are advanced in time. However, note that in
actual numerical simulations $S$ can be very ill-conditioned.  
Now, since $K_{j}=\sum_{i}X_{i}S_{ij}$, the smallest singular value of $K=(K_1,\dots,K_r)$
is equal to that  of $S$. 
Specifically, this is a problem for momentum conservation as
many problems start with zero or very small momentum. This then changes
over time as the algorithm selects appropriate basis functions which
carry a non-zero momentum. However, since initially the contribution
of these functions to $K$ (contrary to $X$) is very small, the coefficients
in the correction have to become large. This implies that the correction
overall becomes quite large. Choosing $X_{j}^{0}$ instead of $K_{j}$,
as we have done here, solves this issue. The situation is analogous
for $V_{i}^{0}$ and $L_{i}$.
\end{rem}

\begin{rem}
\label{rem:Rcorrection}Let us now consider a correction $R_{i}(v)$
for equation (\ref{eq:evol-K})
\[
\partial_{t}L_{i}=\langle X_{i}^{1},F(f)\rangle_{x}+R_{i}(v).
\]
This correction is more general than the ansatz we made in equation
(\ref{eq:ansatz}). As before, our goal is to determine the smallest
$R_{i}(v)$ such that the local conservation laws are satisfied. Since
$R_{i}(v)$ has more degrees of freedom, after the space discretization
has been performed, in principle, a smaller correction could be obtained.
Thus, this seems like a promising approach. In this remark we will
restrict ourselves, for simplicity, only to the continuity equation.
To obey the continuity equation the correction has to satisfy
\[
\int R_{i}(v)\,\mathrm{d}v=-\langle X_{i}^{1},F(f)\rangle_{xv}-\langle X_{i}^{1},\nabla\cdot(\rho u)\rangle_{x}.
\]
Minimizing the correction in the $L^{2}$ norm immediately yields
\[
R_{i}=\frac{-1}{|\Omega_{v}|}\left[\langle X_{i}^{1},F(f)\rangle_{xv}+\langle X_{i}^{1},\nabla\cdot(\rho u)\rangle_{x}\right],
\]
where $|\Omega_{v}|$ is the volume of the domain in the $v$-direction.
Note, in particular, that $R_{i}$ is independent of $v$. Thus, the
correction equally distributes the defect in velocity space. In the
case of the Vlasov equation, however, the density function $f$ is
expected to decay to zero for large velocities. On the other hand,
the described correction would introduce non-zero densities for large
velocities, which is clearly an unphysical artefact. The correction
considered in this paper, i.e.~equation (\ref{eq:ansatz}), only
allows linear combination of $V_{j}^{0}$. This avoids the problems
stated above as the $V_{j}^{0}$ are already used to represent the
numerical solution and thus decay to zero. In fact, any property of
the $V_{j}^{0}$ that is invariant under taking linear combinations,
is preserved by our approach.
\end{rem}

\section{Global conservation\label{subsec:Global-conservation}}

The algorithm developed above satisfies a projected version of the local
conservation law. For mass conservation this is stated as
\[
P_{\overline{X}}\left(\partial_{t}\rho+\nabla\cdot(\rho u)\right)=0.
\]
However, contrary to the continuous formulation, conservation of mass
cannot be deduced from this expression by simply integrating in $x.$
In fact, conservation of mass, in general, is violated for the scheme
described in the previous section. The situation for momentum is similar.

Since we have an underdetermined system of equations it is, in principle,
possible to add an equation that enforces global conservation of mass
and momentum. This has to be done for each step in the splitting algorithm.

\textbf{Step 1: }We impose
\[
0=\partial_{t}\int\rho\,\mathrm{d}x=\sum_{ij}\kappa_{i}\lambda_{ij}\alpha_{j}+\sum_{j}\langle V_{j}^{0},F(f)\rangle_{xv},
\]
where $\alpha_{j}=\int V_{j}^{0}\,\mathrm{d}v$ and $\kappa_{i}=\int X_{i}^{0}\,\mathrm{d}x$,
and
\[
0=\partial_{t}\int\rho u\,\mathrm{d}x=\sum_{ij}\kappa_{i}\lambda_{ij}\beta_{j}+\sum_{j}\langle V_{j}^{0},F(f)\rangle_{xv}\beta_{j},
\]
where $\beta_{j}=\int vV_{j}^{0}\,\mathrm{d}v$. This adds $1+d$
linear equations to the $2(1+d)r$ linear equations (\ref{eq:linear-system})
required for the local conservation laws. Note that in contrast to
these equations all the $\lambda_{ij}$ are coupled to each other. Thus,
we have to solve a single system of size $(1+d)(2r+1)\times r^{2}$.
We will discuss the computational ramifications later in this section. 

\textbf{Step 2: }We impose
\[
0=-\partial_{t}\int\rho\,\mathrm{d}x=\sum_{ij}\kappa_{i}\lambda_{ij}\alpha_{j}+\sum_{ij}\kappa_{i}\langle X_{i}^{1}V_{j}^{0},F(f)\rangle_{xv}\alpha_{j},
\]
where $\alpha_{j}=\int V_{j}^{0}\,\mathrm{d}v$ and $\kappa_{i}=\int X_{i}^{1}\,\mathrm{d}x$,
and
\[
0=-\partial_{t}\int\rho u\,\mathrm{d}x=\sum_{ij}\gamma_{i}\lambda_{ij}\beta_{j}+\sum_{ij}\kappa_{i}\langle X_{i}^{1}V_{j}^{0},F(f)\rangle_{xv}\beta_{j},
\]
where $\beta_{j}=\int vV_{j}^{0}\,\mathrm{d}v$. 

\textbf{Step 3: }We impose
\[
0=\partial_{t}\int\rho\,\mathrm{d}x=\sum_{ij}\kappa_{i}\lambda_{ij}\alpha_{j}+\sum_{i}\kappa_{i}\langle X_{i}^{1},F(f)\rangle_{xv},
\]
where $\alpha_{j}=\int V_{j}^{0}\,\mathrm{d}v$ with $\kappa_{i}=\int X_{i}^{1}\,\mathrm{d}x$,
and
\[
0=\partial_{t}\int\rho u\,\mathrm{d}x=\sum_{ij}\kappa_{i}\lambda_{ij}\beta_{j}+\sum_{i}\kappa_{i}\langle vX_{i}^{1},F(f)\rangle_{xv},
\]
where $\beta_{j}=\int vV_{j}^{0}\,\mathrm{d}v$. 

The problem with this approach is that there is no guarantee that
the resulting linear system even has a solution. This is most easily
demonstrated by considering step 3 in our algorithm. In this case
$b_{i}=0$ (see section \ref{subsec:local-conservation}). 
Now, let us consider the
rank $2$ function on the domain $[0,2\pi]\times\mathbb{R}$ given
by
\[
X_{1}^{1}(x)=\frac{2}{\sqrt{3\pi}}\cos^{2}x,\qquad X_{2}^{1}(x)=\frac{1}{\sqrt{\pi}}\sin(2x),\qquad V_{1}^{0}(v)=\frac{\mathrm{e}^{-v^{2}}}{(\pi/2)^{1/4}},\qquad V_{2}^{0}(v)=\frac{2v\mathrm{e}^{-v^{2}}}{(\pi/2)^{1/4}}.
\]
This gives $\alpha=(\sqrt[4]{2\pi},0)$ and $\kappa=(2\sqrt{\pi/3},0)$.
Thus,
\[
\lambda_{11}=0.
\]
Since $\beta_{1}=0$, we have 
\begin{align*}
\langle X_{1}^{1},F(f)\rangle_{xv} & =\langle X_{1}^{1},\nabla X_{2}^{1}\rangle\cdot\beta_{2}\\
 & \propto\int\cos^{2}(x)\cos(2x)\,\mathrm{d}x\\
 & \neq0.
\end{align*}
This is in contradiction to the condition of global mass conservation. Thus, it is not possible to both satisfy
the continuity equation and obtain global conservation of mass. We
have only considered conservation of mass here, but the same behavior
is observed for the momentum as well. We have the following options:

\textbf{Local:} We enforce only the local conservation laws, while
minimizing the Euclidean norm of the correction.

\textbf{Global:} We enforce only the global conservation laws, while
minimizing the Euclidean norm of the correction.

\textbf{Combined:} We try to find the best approximation to both the
local conservation laws and the global conservation of mass and momentum.
This results in a linear least square problem for the correction.
The different equations can be weighted to either focus
on the local conservation laws or the global conservation of mass
and momentum.

All of these configurations will be considered in section \ref{sec:numerical-results}.
However, before proceeding let us discuss the computational cost of
the combined approach. We have to compute an underdetermined (but
incompatible) linear least square problem with $r^{2}$ unknowns and
$(1+d)(2r+1)$ data. This problem can be solved by computing the Moore\textendash Penrose
pseudo-inverse which requires a QR decomposition of $A^{T}$. Thus,
it requires at most $\mathcal{O}\left(r^{4}\right)$ arithmetic operations
which is typically small compared to the cost of the low-rank algorithm
itself. 

\section{Efficient implementation\label{sec:Efficient-implementation}}

In the proposed algorithm correction terms are added to the three
evolution equations. This implies that our correction is a continuous
function of time for the respective subflows. However, in order to
increase performance it is often of interest to use a specifically
tailored numerical method for solving these subflows. For example,
methods based on fast Fourier techniques (FFT) and semi-Lagrangian
schemes have been proposed in \cite{einkemmer2018low}. To employ these
algorithms while still maintaining the conservation laws for mass
and momentum is not necessarily straightforward. Thus, we will now
introduce a procedure that allows us to apply our correction independent
of the specific time integration strategy that is chosen for solving
the evolution equations (\ref{eq:evolution-corr-K}), (\ref{eq:evolution-corr-S}),
(\ref{eq:evolution-corr-L}). The approach outlined here is similar
to the projection schemes described in \cite{dedner2002}.

We start with the evolution equation for $K_{j}$ which is given as
follows
\begin{align*}
\partial_{t}K_{j} & =\langle V_{j}^{0},F(f)\rangle_{v}+\sum_{k}\lambda_{kj}X_{k}^{0}, \qquad\text{for}\quad f=\sum_l K_l V_l^0.
\end{align*}
Now, we split this equation into
\begin{equation}
\partial_{t}K_{j}=\langle V_{j}^{0},F(f)\rangle_{v}\label{eq:evol-pure-K}
\end{equation}
and
\begin{equation}
\partial_{t}K_{j}=\sum_{k}\lambda_{kj}X_{k}^{0}.\label{eq:evol-corr-step}
\end{equation}
Equation (\ref{eq:evol-pure-K}) is identical to what has to be solved
in case of the original low-rank algorithm described in section \ref{subsec:local-conservation}
(i.e.~the algorithm without correction). Thus, starting from an appropriate
initial value $K_{j}^{0}$ we compute an approximation at time $\tau$,
where $\tau$ is the time step size. This value is henceforth denoted
by $K_{j}^{\star}$. Now, instead of solving equation (\ref{eq:evol-corr-step})
we consider the following approximation
\[
\frac{K_{j}^{1}-K_{j}^{\star}}{\tau}=\sum_{k}\lambda_{kj}X_{k}^{0}.
\]
It remains to derive the conditions on $\lambda_{kj}$ under which the (discretized)
conservation laws are satisfied. We have
\begin{align*}
\rho^{1}-\rho^{0}+\tau\nabla\cdot(\rho^{0}u^{0}) & =\sum_{j}(K_{j}^{1}-K_{j}^{0})\alpha_{j}+\tau\sum_{j}(\nabla K_{j}^{0})\cdot\beta_{j}\\
 & =\sum_{j}(K_{j}^{\star}-K_{j}^{0}+\tau\sum_{k}\lambda_{kj}X_{k}^{0})\alpha_{j}+\tau\sum_{j}(\nabla K_{j}^{0})\cdot\beta_{j},
\end{align*}
where $\alpha_{j}=\int V_{j}^{0}\,\mathrm{dv}$ and $\beta_{j}=\int vV_{j}^{0}\,\mathrm{d}v$.
Now, we apply the projection onto $\overline{X}^0$ to obtain
\[
0=P_{\overline{X}^0}(\rho^{1}-\rho^{0}+\tau\nabla\cdot(\rho^{0}u^{0}))=\sum_{k}X_{k}^{0}\left[\tau\sum_{k}\lambda_{kj}\alpha_{j}+\sum_{j}\langle X_{k}^{0},K_{j}^{\star}-K_{j}^{0}\rangle_{x}\alpha_{j}+\tau\sum_{j}\langle X_{k}^{0},\nabla K_{j}^{0}\rangle_{x}\cdot\beta_{j}\right].
\]
This is the analogue to equation (\ref{eq:continuity-K}). 

For the momentum balance equation we have
\begin{align*}
\rho^{1}u^{1}-\rho^{0}u^{0}+\tau\nabla\cdot(\rho^{0}u^{0}\otimes u^{0})+\tau E^{0}\rho^{0} & =\sum_{j}(K_{j}^{1}-K_{j}^{0})\beta_{j}+\tau\sum_{j}(\nabla K_{j}^{0})\cdot\gamma_{j}+\tau E^{0}\rho^{0}\\
 & =\sum_{j}(K_{j}^{\star}-K_{j}^{0}+\tau\sum_{k}\lambda_{kj}X_{k}^{0})\beta_{j}+\tau\sum_{j}(\nabla K_{j}^{0})\cdot\gamma_{j}+\tau E^{0}\rho^{0},
\end{align*}
where $\gamma_{j}=\int(v\otimes v)V_{j}^{0}\,\mathrm{d}v$ and we
have used $E^{0}$ to denote the electric field at the beginning of
the time step. Applying the projection onto $\overline{X}^0$ we obtain
\begin{align*}
0 & =P_{\overline{X}^0}(\rho^{1}u^1-\rho^{0}u^0+\tau\nabla\cdot(\rho^{0}u^{0}\otimes u^{0})+\tau E^{0}\rho^{0})\\
 & =\sum_{k}X_{k}^{0}\left[\tau\sum_{j}\lambda_{kj}\beta_{j}+\sum_{j}\langle X_{k}^{0},K_{j}^{\star}-K_{j}^{0}\rangle_{x}\beta_{j}+\tau\sum_{j}\langle X_{k}^{0},\nabla K_{j}^{0}\rangle_{x}\cdot\gamma_{j}+\tau\langle X_{k}^{0},E^{0}\rho^{0}\rangle_{x}\right]
\end{align*}
which is the analogue to equation (\ref{eq:momentumbal-K}). 

In fact, these equations are precisely in the form of (\ref{eq:linear-system}).
Only the right-hand side
\begin{align*}
b_{i} & =-\sum_{j}\biggl\langle X_{k}^{0},\frac{K_{j}^{\star}-K_{j}^{0}}{\tau}\biggr\rangle_{x}\alpha_{j}-\sum_{j}\langle X_{k}^{0},\nabla K_{j}^{0}\rangle_{x}\cdot\beta_{j}\\
d_{i} & =-\sum_{j}\biggl\langle X_{k}^{0},\frac{K_{j}^{\star}-K_{j}^{0}}{\tau}\biggr\rangle_{x}\beta_{j}-\sum_{j}\langle X_{k}^{0},\nabla K_{j}^{0}\rangle_{x}\cdot\gamma_{j}-\sum_{j}\langle X_{k}^{0},E^{0}K_{j}^{0}\rangle_{x}\alpha_{j}
\end{align*}
 is modified. Thus, there is no additional difficulty in implementing
this approach.

A similar procedure can be applied to step 2 and 3 of the splitting
algorithm. For step 2 we obtain equation (\ref{eq:linear-system})
with
\begin{align*}
b_{i} & =\sum_{j}\frac{S_{ij}^{\star}-S_{ij}^{1}}{\tau}\alpha_{j}-\sum_{j}\langle X_{i}^{1},\nabla X_{j}^{1}\rangle_{x}S_{ij}^{1}\beta_{j}\\
d_{i} & =\sum_{j}\frac{S_{ij}^{\star}-S_{ij}^{1}}{\tau}\beta_{j}-\sum_{j}\langle X_{j}^{1},\nabla X_{i}^{1}\rangle_{x}S_{ij}^{1}\gamma_{j}-\sum_{j}\langle E^{0},X_{j}^{1}\rangle_{x}S_{ij}^{1}\alpha_{j}
\end{align*}

and for step 3 we obtain equation (\ref{eq:linear-system}) with (denoting the integral $\la g \ra_v = \int g \,\mathrm{d}v$)
\begin{align*}
b_{i} & =-\biggl\langle\frac{L_{i}^{\star}-L_{i}^{0}}{\tau}\biggr\rangle_{v}-\sum_{j}\langle X_{i}^{1},\nabla X_{j}^{1}\rangle_{x}\cdot\langle vL_{j}^{0}\rangle_{v}\\
d_{i} & =-\biggl\langle\frac{L_{i}^{\star}-L_{i}^{0}}{\tau}\biggr\rangle_{v}-\sum_{j}\langle X_{i}^{1},\nabla X_{j}^{1}\rangle_{x}\cdot\langle(v\otimes v)L_{j}^{0}\rangle_{v}-\sum_{j}\langle X_{i}^{1},EX_{j}^{1}\rangle_{x}\langle L_{j}^{0}\rangle_{v}.
\end{align*}

Thus, we are able to apply the procedure introduced in section \ref{subsec:local-conservation},
independent of the specific numerical discretization. This has the
added benefit that the correction and the associated coefficients
only need to be computed once for each step of the splitting algorithm.
The only downside here is that we have traded the continuous version
of the conservation laws for a discretized version.

\section{Numerical results\label{sec:numerical-results}}

In this section we will present numerical results for a two-stream
instability. Specifically, we consider the domain $[0,10\pi]\times[-9,9]$
and impose the initial value
\[
f_{0}(x,v)=\frac{1}{2\sqrt{2\pi}}\left(\mathrm{e}^{-(v-v_{0})^{2}/2}+\mathrm{e}^{-(v+v_{0})^{2}/2}\right)(1+\alpha\cos(kx)),
\]
where $\alpha=10^{-3}$, $k=\tfrac{1}{5}$, and $v_{0}=2.4$. Periodic
boundary conditions are used in both the $x$- and the $v$-direction.
This setup models two beams propagating in opposite directions and is an unstable equilibrium. Small perturbations in the initial particle-density
function eventually force the electric energy to increase exponentially.
This is called the linear regime. At some later time saturation sets
in (the nonlinear regime). This phase is characterized by nearly constant
electric energy and significant filamentation of the phase space.
This test problem has been considered in \cite{Ehrlacher2017,Kormanna} and \cite{einkemmer2018low}
in the context of low-rank approximations. It has been established
there that low-rank approximations of relatively small rank are sufficient
in order to resolve the linear regime. However, once saturation sets
in, the reference solution (computed using a full grid simulation)
shows only small oscillations in the electric field. For the low-rank
approximation, however, oscillations with significant amplitude can
be observed. Since filamentation makes it very difficult to efficiently
resolve the small structures in this regime (the $L^{\infty}$
error will be large for any numerical method), we consider it a good
test example for the conservative method developed in this work.

\begin{figure}
\begin{centering}
\includegraphics[width=14cm]{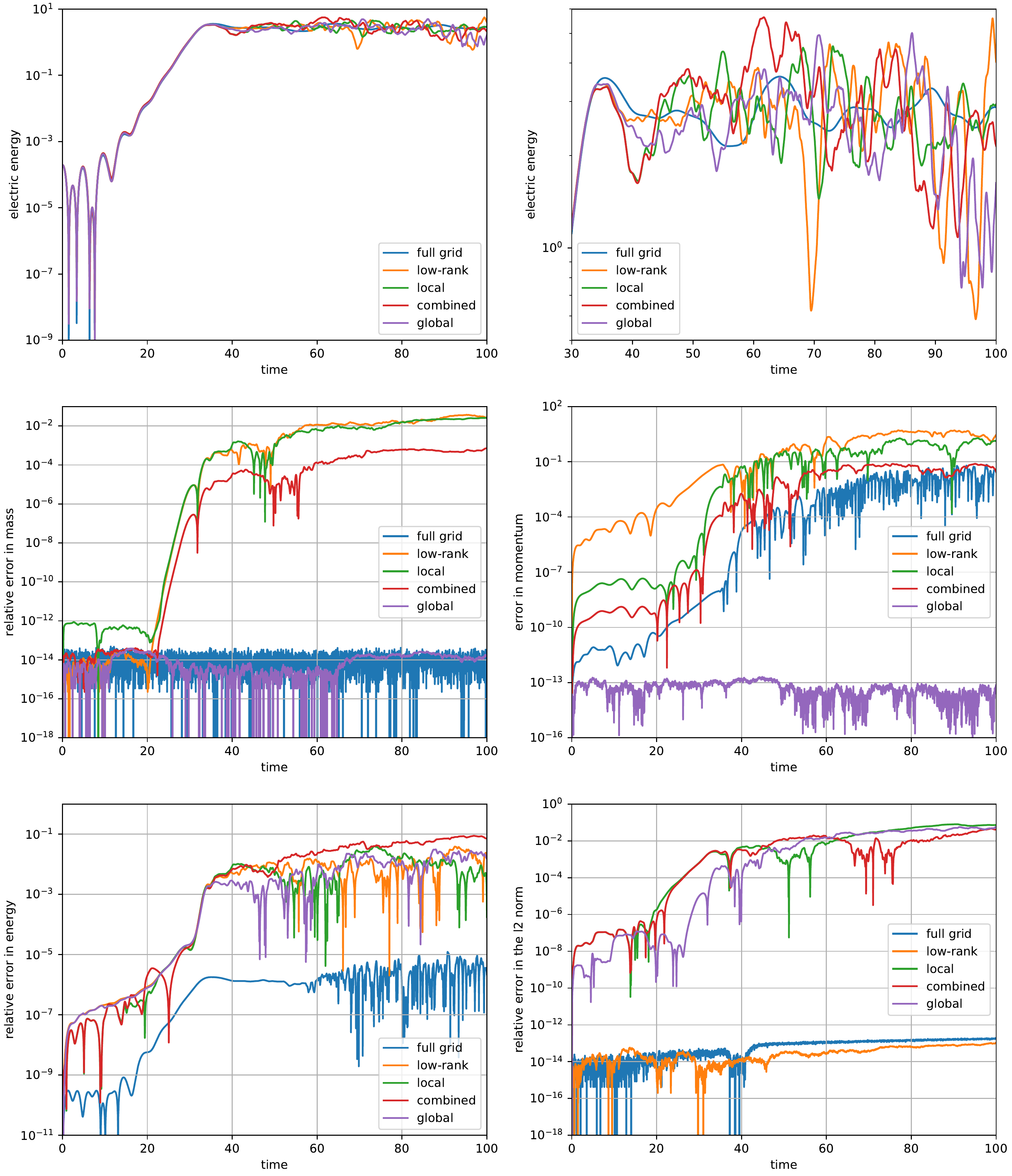}
\par\end{centering}
\caption{Numerical simulations of the two-dimensional two-stream instability
with rank $r=10$ are shown. The Strang splitting algorithm with a
time step size $\tau=0.025$ is employed. In both the $x$ and $v$-directions
$128$ grid points are used. As a comparison, a direct Eulerian simulation
(based on a spectral method) is also shown. \label{fig:ts-r10}}
\end{figure}

In Figure \ref{fig:ts-r10} numerical simulation of the two-stream
instability for rank $r=10$ are shown for the algorithm without correction
(labeled low-rank), the correction that exactly satisfies the local projected 
continuity equations described in section \ref{subsec:local-conservation}
(labeled local), the algorithm of section \ref{subsec:Global-conservation}
that combines both local and global corrections (labeled combined),
and the algorithm that conserves mass and momentum exactly but does
not satisfy the local continuity equations (labeled global). In addition,
the full grid simulation is shown (labeled full grid). We observe
that all methods show excellent agreement in the linear regime. In
the nonlinear regime the local correction shows the best performance
(the least amount of oscillations). The performance of the combined
approach is also significantly better compared to the uncorrected
algorithm and the global correction. The uncorrected algorithm clearly
performs worst.

Figure \ref{fig:ts-r10} also shows the error in mass, momentum, energy,
and the $L^{2}$ norm. We see that although the local correction results
in a significant improvement with respect to the qualitative behavior
of the electric field, the errors in mass and momentum are still comparable
to the uncorrected algorithm. As has been discussed in section \ref{subsec:Global-conservation},
in general, satisfying both the local continuity equations and the
global invariants is not possible. We clearly see this in the numerical
simulation. Nevertheless, the combined approach results in a significant
reduction in the error in mass and momentum (by approximately two
orders of magnitude).

\begin{figure}
\begin{centering}
\includegraphics[width=14cm]{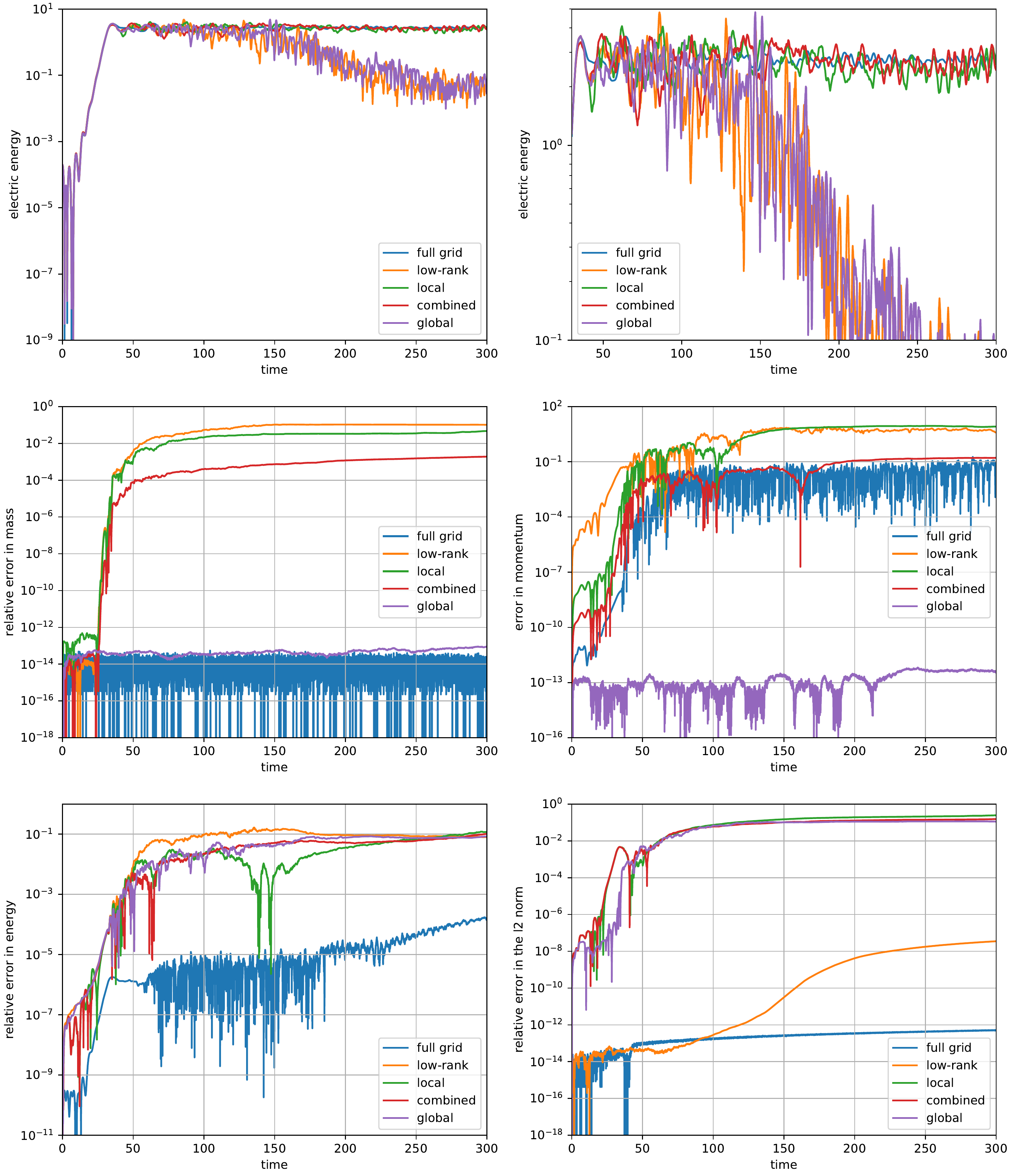}
\par\end{centering}
\caption{Numerical simulations of the two-dimensional two-stream instability
with rank $r=15$ are shown. The Strang splitting algorithm with a
time step size $\tau=0.025$ is employed. In both the $x$ and $v$-directions
$128$ grid points are used. As a comparison, a direct Eulerian simulation
(based on a spectral method) is also shown. \label{fig:ts-r15}}
\end{figure}

Now, we increase the rank to $r=15$ and consider a longer time interval
(up to $t=300$). The numerical results are shown in Figure \ref{fig:ts-r15}.
It can be observed very clearly that the uncorrected algorithm as
well as the global correction result in qualitatively wrong results
(the electric energy decreases by more than two orders of magnitude).
On the other hand, the local correction and the combined approach keep
the electric energy stable until the final time of the simulation.
With respect to the conservation of the invariants the same conclusion
as above can be drawn.

\begin{figure}
\begin{centering}
\includegraphics[width=14cm]{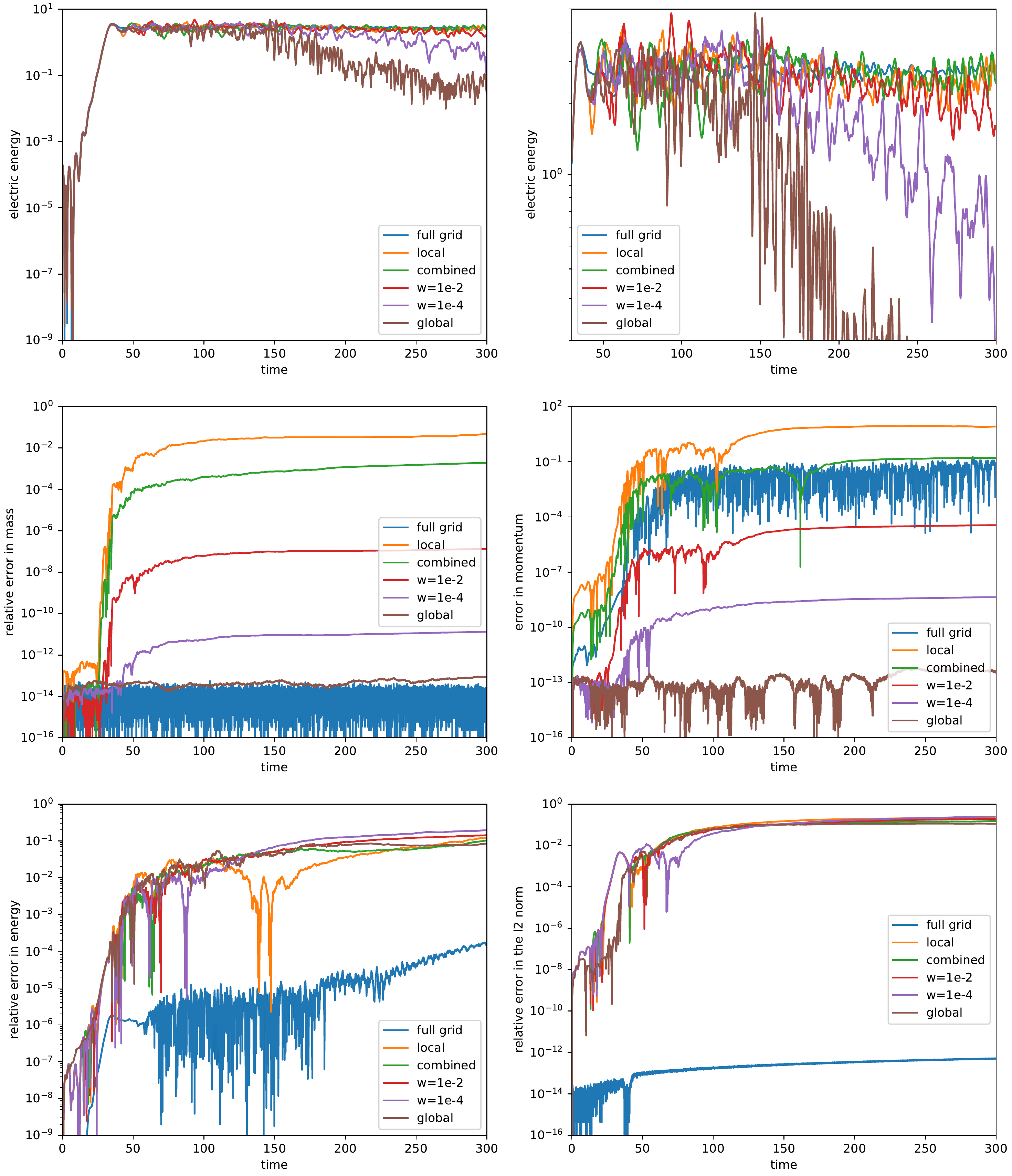}
\par\end{centering}
\caption{Numerical simulations of the two-dimensional two-stream instability
with rank $r=15$ are shown. For the combined approach we also show
numerical results for $w=10^{-2}$ and $w=10^{-4}$. A weight of $w=0$
corresponds to the global correction and a weight of $w=1$ to the
combined correction described in section \ref{subsec:Global-conservation}.
The Strang splitting algorithm with a time step size $\tau=0.025$
is employed. In both the $x$ and $v$-directions $128$ grid points
are used. As a comparison, a direct Eulerian simulation (based on
a spectral method) is also shown.\label{fig:tsw-r10}}
\end{figure}

As has been mentioned in section \ref{subsec:Global-conservation},
the combined approach can be adjusted to either be closer to the local
correction or the global correction. The results in Figure \ref{fig:tsw-r10}
show how we can trade-off the error in mass and momentum and the error
in the local conservation laws. We clearly see that the solution deteriorates
as the error in the conservation laws increases.

%\section*{References}

\bibliographystyle{plain}
\bibliography{vlasov-lrcons}

\end{document}